\documentclass[11pt]{article}

\newtheorem{Theo}{Theorem}[section]
\newtheorem{Lem}[Theo]{Lemma}
\newtheorem{Pro}[Theo]{Proposition}
\newtheorem{rem}[Theo]{Remark}

\newtheorem{Ex}[Theo]{Example}
\newtheorem{Def}[Theo]{Definition}

\begin{document}

\sf
\title{HOPF QUIVERS}
\author{Claude CIBILS - Marc ROSSO}
\date{ }

\maketitle

\begin{abstract}
We classify graded Hopf algebras structures over path coalgebras,
that is over free pointed coalgebras, using Hopf quivers which are
analogous to Cayley graphs.
\end{abstract}

\small
\noindent 2000 Mathematics Subject Classification : 57T05
16W30

\noindent Keywords : Hopf, algebra, coalgebra, quiver,
path, pointed.

\normalsize

\section{\sf Introduction}

In this paper we provide the classification of path
coalgebras which admit a graded Hopf algebra structure; we
construct pointed free Hopf algebras in an exhaustive way.
Paths of a quiver provide a natural basis of such algebras
and the comultiplication of a path is the sum of all the
splits of the path. We show that the multiplication of two
paths is described through {\em thin splits} of the paths,
a group structure on the vertices and an action of this
group on the arrows by permutations on the left and
permutations and linear action on the right.

In a previous paper \cite{ciro} we have presented parts of
these results in a dual version since we considered
associative path algebras provided by finite quivers. In
\cite{ciro}, we have described basic Hopf algebras (simple
modules are one-dimensional) instead of pointed Hopf
algebras. However this approach appears to be somehow
unusual. Note that Radford \cite{ra} describes the
structure of finite dimensional simple-pointed Hopf
algebras generated by pairs of a group-like $a$ and a
$(1,a)$-primitive, when the field is algebraically closed.

In the present paper we obtain a complete and more suitable
version of the classification which allows infinite
quivers. Moreover the description involves formulas for the
product besides the canonical formulas for the coproduct.
This makes explicit the quantum shuffle product \cite{ro}
making use of natural elements in a Hopf bimodule having a
simple "geometrical" interpretation in the quiver sense,
rather than working systematically with right or left
coinvariants.

 The first author lectured at MSRI \cite{cimsri} on the
present approach in November 1999. The books of C. Kassel
\cite{ka}, S. Montgomery \cite{mo}, and M.E. Sweedler
\cite{sw} provide the necessary background for the Hopf
algebra side of this paper, while the book of M. Auslander,
I. Reiten and S. Smal\o \ \cite{auresm} gives precise
insight to the path algebra approach.

\section{\sf Path Coalgebras}

A quiver $Q$ is an oriented graph given by two sets $Q_{0}$
and $Q_{1}$ of vertices and arrows and two maps $s, t :
Q_{1} \to Q_{0}$ providing each arrow with its source and
terminal vertices. Infinite sets $Q_{0}$ and $Q_{1}$ are
allowed. A path $\alpha$ is a finite sequence of
concatenated arrows $\alpha = a_{n}\cdots a_{1}$ which
means that $t(a_{i}) = s(a_{i+1})$ for $i = 1,\cdots,n-1$.
We set $s(\alpha) = s(a_{1})$ and $t(\alpha) = t(a_{n})$.
Moreover, a vertex $u$ coincides with its source and
terminal vertices. The length of a path is the length of
its arrow sequence; vertices are zero-length paths.

If one prefers, a quiver is the same structure than a free
category with a given set of free generators : $Q_{0}$ is
the set of objects, the set of paths is the set of
morphisms and $Q_{1}$ is the free set of generators. A
functor $F$ from this category to another one is completely
determined by the image of the objects and a coherent
choice of images of the arrows, meaning that for each arrow
$a$ the morphism $F(a)$ is from $F(s(a))$ to $F(t(a))$.

Let now $k$ be a field.

\begin{Def}
The coalgebra of paths of a quiver $Q$ is the linearisation
$kQ$ of the set of paths, equipped with a comultiplication
$\Delta$ and a counit $\varepsilon$ as follows. For a path
$\alpha = a_{n} \cdots a_{1}$ we set $$
\begin{array}{lll}
\Delta \alpha & =& \sum^n_{i=0} a_{n} \cdots a_{i+1}
\otimes a_{i} \cdots a_{1} \\ & = & \alpha \otimes
s(\alpha) + \sum_{i=1}^{n-1} a_{n} \cdots a_{i+1} \otimes
a_{i} \cdots a_{1} + t(\alpha) \otimes \alpha \\ & = &
\sum_{\alpha = \alpha_{(2)} \alpha_{(1)}} \alpha_{(2)}
\otimes \alpha_{(1)}
\end{array}
$$

where $\{\alpha_{(2)} \alpha_{(1)}\}$ is the set of splits of
$\alpha$. Moreover $\varepsilon(\alpha) = 0$ if $\alpha$ is a
positive length path and $\varepsilon(u) = 1$ if $u$ is a vertex.
\end{Def}

There is no difficulty to verify that $kQ$ is indeed a
coalgebra. Actually this coalgebra structure is precisely
the dual version of the well known path algebra structure
on finite quivers as it is widely used. It has been
considered previously by  W. Chin and S. Montgomery
\cite{chmo}, see also \cite{gr,ta,ni} and \cite{tooyzh}.

We also note that a path coalgebra is an instance of a
cotensor coalgebra. Indeed consider the linearisation
$kQ_{0}$ of the set of vertices $Q_{0}$ with its natural
coalgebra structure, namely each vertex is a group-like.
The vector space $kQ_{1}$ is a $kQ_{0}$-bicomodule
\emph{via} $$\delta_{L}(a) = t(a) \otimes a \hskip1cm
\delta_{R}(a) = a \otimes s(a).$$

The cotensor product of $kQ_{1}$ with itself is the kernel of
$\delta_{R} \otimes 1 - 1 \otimes \delta_{L}$.

We denote $Q_{n}$ the set of paths of length $n$. One
immediately notice that the cotensor square of $kQ_{1}$ is
$kQ_{2}$, where the later is considered with its natural
$kQ_{0}$-bicomodule structure. In other words, the cotensor
coalgebra over the trivial coalgebra $kQ_{0}$ of the
bicomodule $kQ_{1}$ is the path coalgebra we have described
above.

\begin{Def}
Let $x$ and $y$ be vertices. The \emph{$(y,x)$-isotypic
component} of a $kQ_0$-bicomodule $Z$ is
 $$ ^yZ^{x} =
\{ z\in Z\ \vert \ \delta_{L}(z) = y \otimes z \quad {\rm
and} \quad \delta_{R}(z) =z \otimes x\}.$$ In case of
$kQ_n$ it corresponds to the vector space of $n$-paths from
$x$ to $y$.
\end{Def}

We record now a useful consequence, namely the universal
property enjoyed by $kQ$. A coalgebra map $\psi : X \to kQ$
is given by a sequence of maps $\psi_n : X\to kQ_n$, and is
uniquely determined by two maps $\psi_{0} : X \to kQ_{0}$ a
coalgebra map and $\psi_{1} : X \to kQ_{1}$ a bicomodule
map, where $X$ is to be considered as a $kQ_{0}$-bicomodule
via $\psi_{0}$. For instance, the map $\psi_{2} : X \to
kQ_{2}$ is provided by the composition $$X
\stackrel{\Delta}{\longrightarrow} X \otimes X
\stackrel{\psi_{1} \otimes \psi_{1}}{\longrightarrow}
kQ_{1} \otimes kQ_{1}$$ which image is indeed contained in
he cotensor square $kQ_{2}$. More generaly $\psi_{n}
=\psi_{1}^{\otimes n}\circ \Delta^{(n)}$ where
$\Delta^{(n)}$ denotes $n$ successive applications of the
comultiplication, a well defined map since $\Delta$ is
coassociative.

\section{\sf Hopf quivers}

Hopf quivers are defined bellow, they are precisely the
quivers such that the path coalgebra can be endowed with a
graded Hopf algebra structure. Those quivers are similar to
Cayley graphs, which have set of vertices given by the
elements of a group and arrows corresponding to
multiplication by elements of a chosen system of
generators.

A {\em ramification data} $r$ of a group $G$ is a positive
central element of the group ring of $G$, {\it i.e.} $r =
\sum_{c\in \mathcal{C}} r_{C} C$ is a formal sum of
conjugacy classes of $G$ with positive integer
coefficients.

\begin{Def}
Let $G$ be a group and $r$ be a ramification data. The
corresponding Hopf quiver has set of vertices the elements
of $G$ and has $r_{C}$ arrows from $x$ to $cx$ for each $x
\in G$ and $c \in C$.
\end{Def}

A Hopf quiver is connected if and only if the union of the
conjugacy classes affording non zero coefficients in the
ramification data generates $G$.

\begin{Ex}\sf
Let $G$ be a cyclic group which can be infinite, let $K$ be
a generator of $G$ and $r = K^{2}$. The Hopf quiver of $G$
is connected if $|G|$ is finite and odd, otherwise it has
two connected components. In case $G$ is finite the
connected components are crowns.
\end{Ex}

\begin{Theo}
\label{classe}
Let $Q$ be a quiver. The path coalgebra $kQ$
admits a graded Hopf algebra structure if and only if $Q$
is the Hopf quiver of some group with respect to a
ramification data.
\end{Theo}

\noindent \textbf{Proof} : If the path coalgebra $kQ$ is a
graded Hopf algebra, the set $Q_{0}$ of vertices is the set
of group-likes, hence $Q_{0}$ is a group and the
sub-coalgebra $kQ_{0}$ turns into a sub-Hopf algebra,
namely the group algebra of $Q_{0}$. Moreover the arrows
provide a sub-vector space $k Q_{1}$ which is
simultaneously a $kQ_{0}$-bimodule and a
$kQ_{0}$-bicomodule with structure maps compatible, since
$\delta_{L}$ and $\delta_{R}$ are $kQ_{0}$-bimodule maps.

Such structures are called {\em Hopf bimodules}; we recall
in the next section their classification. It turns out that
the dimension of an isotypic component
 $$ ^y(kQ_{1})^{x} =
\{ \alpha \in kQ_{1}\ \vert \ \delta_{L}(\alpha) = y \otimes
\alpha \quad {\rm and} \quad \delta_{R}(\alpha) =\alpha \otimes
x\} $$ is constant when $yx^{-1}$ remains in the same conjugacy
class. In other words the number of arrows from $x$ to $y$ only
depends on the conjugacy class of $yx^{-1}$. Consequently the
quiver is the Hopf quiver of the ramification data which
coefficients are the dimensions of the isotypic components.

Conversely, let $Q$ be the Hopf quiver of a group $G$ with respect
to a ramification data $r$. We will show in the next section that
there always exist a $kG$-Hopf bimodule $B$ such that the
dimensions of the isotypic components of $B$ are prescribed by
$r$, namely such that $$ \dim_{k}\ \ ^y\! B^x =r_{\langle
yx^{-1}\rangle}.   $$

Actually there are several Hopf bimodules affording the
same ramification data. We chose one of them and we provide
the vector space $kQ_{1}$ with the transported Hopf
bimodule structure : chose a basis of $B$ which respects
the isotypic components and identify the vector spaces $B$
and $kQ_{1}$ using the set of arrows as a basis of
$kQ_{1}$. At this point we have a group structure on the
vertices of $Q$ and a compatible multiplication of vertices
by arrows. Using the quoted universal property of path
coalgebras this data uniquely determines an associative
coalgebra morphism and $kQ$ becomes a pointed bialgebra.
The existence of the antipode is granted by Takeuchi's
results \cite{ta} p. 572 since the group-like elements of
$kQ$ are invertible.

Next we display the formula for the product of two arrows
which is obtained directly from the description we made at
the end of the preceding section of the map in degree 2.
Let $a$ and $b$ be arrows of the quiver. Then $$ a.b =
[t(a).b][a.s(b)] + [a.t(b)][s(a).b]. $$

Note  for instance that $t(a).b$ denotes the left action of the
group element $t(a)$ on the arrow $b$. We record that the two
homogeneous  terms of $a.b$ starts at the vertex $s(a) s(b)$ and
ends at $t(a)t(b)$. The above formula is a particular case of the
following Theorem which describes the product of two paths.

\begin{Def}
A $p$-thin split of a path $\alpha$ is a sequence $\alpha_{p},
\alpha_{p-1},\cdots,\alpha_{1}$ of vertices and arrows such that
the product $\alpha_{p} \alpha_{p-1} \cdots \alpha_{1}$ (in the
path algebra sense) is $\alpha$.
\end{Def}

\begin{Ex}
{\sf Let $\alpha = cba$ be a length $3$ path. A $7$-thin cut of
$\alpha$ is for instance $(t(c),t(c),c,b,s(b),s(b),a)$.}
\end{Ex}

Clearly $p$-thin splits of a  $n$-path $\alpha$ are in
one-to-one correspondence with $p$-sequences of zeros and
ones such that the number of ones is $n$. We denote
$D^p_{n}$ the set of such sequences. More precisely, if $d
\in D^p_{n}$ and $(a_{n}, \cdots, a_{1})$ is the sequence
of arrows of $\alpha$, the p-thin split $d\alpha =
((d\alpha)_{p}, \cdots, (d\alpha)_{1})$  is determined by
``$(d\alpha)_{i}$ is a vertex if $d_{i} = 0$ and is an
arrow if $d_{i} = 1$". The 7-thin cut of the example above
corresponds to $d = (0,0,1,1,0,0,1)$.

Let $\alpha$ be a n-path and $\beta$ be a m-path. Let $d
\in D_{n}^{n+m}$ and let $\bar{d} \in D_{n}^{n+m}$ be the
complement sequence obtained from $d$ by replacing each $0$
by a $1$ and each $1$ by a $0$. Consider the element
 $$
(\alpha . \beta)_{d} = [(d\alpha)_{m+n} . (\bar{d} \beta)_{m+n}]
\cdots [(d \alpha)_{1} . (\bar{d} \beta)_{1}] $$
 which lies in the
$(m+n)$-cotensor power of $B$ and belongs  to the isotypic
component of type $(t(\alpha)t(\beta),s(\alpha)s(\beta))$.

\begin{rem}
\sf If $d_i=1$, the element $(\bar{d} \beta)_{i}$ is a vertex
which acts on the right on the arrow $(d\alpha)_{i}$ and
$[(d\alpha)_{i} . (\bar{d} \beta)_{i}]$ denotes the result of the
action. Conversely, if $d_{i} = 0$ the above expression is the
result of the left action of the group-vertex element
$(d\alpha)_{i}$ on the arrow $(\bar{d} \alpha)_{i}$.
\end{rem}

\begin{rem} \sf The source of the first term of $(\alpha .
\beta)_{d}$ is $$ s[(d\alpha)_{1} . (\bar{d} \beta)_{1}] =
s((d\alpha)_{1}) s((\bar{d} \beta)_{1}) = s(\alpha) s(\beta) $$
and the sequence of terms of $(\alpha . \beta)_{d}$ is
concatenated.
\end{rem}

\begin{Theo}
Let $G$ be a group, $k$ be a field, $B$ be a $kG$-Hopf
bimodule and let $C$ be the associated cotensor coalgebra
endossed with its Hopf algebra structure (see the proof of
Theorem \ref{classe}). Chose a basis of the isotypic
components of the bicomodule $B$ and identify the coalgebra
$C$ with the path coalgebra of the quiver of $B$. Let
$\alpha$ and $\beta$ be respectively paths of length $n$
and $m$. Then $$ \alpha. \beta = \sum_{d \in D_{n}^{n+m}}
(\alpha . \beta)_{d} $$
\end{Theo}

\begin{rem} $|D_{m}^{n+m}| = ({m+n \atop n})$.
\end{rem}

The proof of the Theorem uses the following result :

\begin{Lem}
Let $\Delta_2$ be the comultiplication of the path coalgebra $kQ
\otimes kQ$. Then $$ \Delta^{(p)}_{2} (\alpha \otimes \beta) =
\sum_{\alpha_{(p)}\cdots \alpha_{(1)} = \alpha \atop \beta_{(p)}
\cdots \beta_{(1)} = \beta} \alpha_{(p)} \otimes \beta_{(p)}
\otimes \cdots \otimes \alpha_{(1)} \otimes \beta_{(1)} $$
\end{Lem}

\noindent \textbf{Proof} : Note that the sum of the formula is
over all the $p$-splits of $\alpha$ and $\beta$. Recall that
$\Delta_{2} = (1 \otimes \tau \otimes 1)(\Delta \otimes \Delta)$
where $\tau$ is the flip map. Hence $$ \Delta_{2}(\alpha \otimes
\beta) = \Sigma \alpha_{(2)} \otimes \beta_{(2)} \otimes
\alpha_{(1)} \otimes \beta_{(1)}. $$ The result follows by
induction.

\noindent \textbf{Proof of the Theorem}: we  use the universal
property of $C$ with respect to $X = C \otimes C$. We already have
in low degrees $$ \varphi_{0} : C \otimes C \to kG \otimes kG \to
kG \subset C $$ $$ \varphi_{1} : C \otimes C \to kG \otimes B + B
\otimes kG \to B \subset C $$ where the first maps for
$\varphi_{0}$ and $\varphi_{1}$ are  respectively the projections
to the 0 and 1 homogeneous component of $C \otimes C$. In order to
describe $\varphi_{n+m}(\alpha \otimes \beta)$, we first use the
Lemma for obtaining  $\Delta^{(n+m)}(\alpha \otimes \beta)$. The
morphism $\varphi_{1}^{\otimes (n+m)}$ only retains pairs of
complementary $(n+m)$-thin splits. The formula follows immediately
from this remark.

\begin{Ex}
\sf If $a$ and $b$ are arrows, the $2$-thin splits of $a$
are $(a,s(a))$ and $(t(a),a)$. The complement thin splits
of $b$ are respectively $(t(b),b)$ and $(b,s(b))$. Then $$
a.b = [a.t(b)] [s(a).b] + [t(a).b][a.s(b)].$$
\end{Ex}

\begin{Ex}
\sf
 Let $a$ be an arrow and $\beta = cb$ be a
$2$-path. The $3$-thin splits of $a$ and the corresponding
complement $3$-thin splits of $\beta$ are $$
\begin{array}{rcccl}
(&t(a),&t(a),&a&) \\ (&c,&b,&s(b)&)\\ \end{array} $$ $$
 \begin{array}{rcccl}
(&t(a),&a,&s(a)&)
\\ (&c,&t(b),&b&)\\
\end{array}
 $$ $$
 \begin{array}{rcccl}
 (&a,&s(a),&s(a)&) \\ (&t(c),&c,&b&)
\end{array}
$$
Then $$
\begin{array}{lll}
a.(cb) & = & [t(a).c][t(a).b][a.s(b)] + \\ & &
[t(a).c][a.t(b)][s(a).b] + \\ & & [a.t(c)][s(a).c][s(a).b].
\end{array}
$$

We record that if $u$ is a vertex and $\alpha = (a_{n}, \cdots,
a_{1})$ is a path, then $$ \begin{array}{l}
  u\cdot
(a_{n} \cdots a_{1}) = [u\cdot a_{n}] \cdots [u \cdot
a_{1}] \\ (a_{n} \cdots a_{1})\cdot  u = [a_{n}\cdot
u]\cdots [a_{1}\cdot u].
\end{array}$$
\end{Ex}

Next we will display a formula for the product of a
sequence of $n$ arrows $A = (a_{n},\cdots,a_{1})$. We
stress that there is no need for the arrows to be
concatenated. Recall that the vertices of the quiver are
group elements which acts on both sides on the vector space
generated by the arrows.

Let $\sigma$ be a permutation of the set $\{1,\cdots,n\}$.
We define a map $p^\sigma_{k}$ which assigns to each arrow
$a_{i}$ either its source or terminus vertex, according to
$i$ being already reached by $\sigma$ acting on
$\left\{1,\cdots,k-1\right\}$:

$$ p^\sigma_{k}(a_{i}) = \left \{ \begin{array}{lll}
s(a_{i}) & {\sf if} & i\not
\in\sigma\left\{1,\cdots,k-1\right.\} \\ t(a_{i}) & {\sf
if} & i\in\sigma\left\{1,\cdots,k-1\right\}.
\end{array} \right.
$$

Given a sequence of arrows $A=(a_{n},\cdots,a_{1})$, a
permutation $\sigma$ and an integer $k \in \{1,\cdots,n\}$
we define an homogeneous element $$ A^\sigma_{k} = \left (
p^{\sigma}_{k}(a_{n}) \cdots
p^{\sigma}_{k}(a_{\sigma(k)+1}) \right ) .\ a_{\sigma(k)}\
. \left (p^{\sigma}_{k} (a_{\sigma(k)-1}) \cdots
p^{\sigma}_{k}(a_{1}) \right ).$$
 Note that for a fixed
$\sigma$ the elements $ A^\sigma_{k}$ are concatenated as
$k$ increases. Finally we set $A^\sigma=A^\sigma_{n} \cdots
A^\sigma_{1}$ where the product is the usual one in the
path algebra sense. $A^\sigma$ is an homogenous element of
degree $n$ from the vertex $s(a_{n})\cdots s(a_{1})$ to the
vertex $t(a_{n}) \cdots t(a_{1})$.

\begin{Pro}
Let $A$ be a sequence of non necessarily concatenated
arrows $(a_{n},\cdots,a_{1})$. Then $$ a_{n} \cdot \cdots
\cdot a_{1} = \sum_{\sigma \in S_{n}} A_{\sigma}. $$
\end{Pro}

The proof follows from the product formula of paths that we
provided above.

\begin{Ex}\sf For $n = 3$, consider $A = (c,b,a)$.
Then $$\begin{array}{lll} c.b.a & = &
[c.tb.ta][sc.b.ta][sc.sb.a] + \\ & &
[tc.b.ta][c.sb.ta][sc.sb.a] + \\ & &
[c.tb.ta][sc.tb.a][sc.b.sa] + \\ & &
[tc.tb.a][tc.b.sa][c.sb.sa] +
\\ & & [tc.tb.a][c.tb.sa][sc.b.sa] + \\ & &
[tc.b.ta][tc.sb.a][c.sb.sa].
\end{array}$$

The first row corresponds to the identity permutation, next we use
the three transpositions and finally the elements of order 3.
\end{Ex}

We provide now  examples of  Hopf algebras obtained through
the procedure of path coalgebras.

\begin{Ex}
\sf Let $G = \{1\}$ be the trivial group and $B = k$ be the
trivial Hopf bimodule. The ramification data is $r = 1$,
and the quiver is a loop $x$. The path coalgebra
corresponds to the polynomial algebra in one variable,
while the path algebra structure corresponds to the usual
product of polynomials. The formula above for the
multiplication of paths provides the product $$ X^n.X^m =
\left ( {n+m \atop m} \right ) X^{n+m} $$
\end{Ex}

\begin{Ex} \sf
Let $G$ be a cyclic group, $K$ be a generator, and consider
the $kG$ Hopf bimodule defined as follows ($q$
 is a non-zero element if $G$ is infinite or a root of
unity if $G$ is finite) : the isotypic components
${}^{K^{i+1}} B^{K^i}$ are one-dimensional with basis
$E_{i}$. Other components ${}^{K^j} B^{K^i}$ are zero. The
left action is by translation $$ K.E_{i} = E_{i+1} $$ while
the right action is by translation and multiplication by q
$$ E_{i}.K = q E_{i+1}. $$ The corresponding Hopf quiver is
a crown or a quiver of type $A^\infty_{\infty}$. As for any
path coalgebra a natural basis is provided by the set of
paths. In this example the basis consists of concatenated
sequences of arrows $E_{i}$ and 0-length paths (i.e.
vertices).
\end{Ex}

The Gauss binomial coefficient $\left ( {n \atop i} \right
)_{q}$ is defined using the algebra $k\{x,y\} / \langle yx
- q xy \rangle$, we have $ (x+y)^n = \sum^n_{0} \left ( {n
\atop i} \right )_{q} x^i y^{n-i}$. The recursive formula
$\left ( {n+1 \atop i} \right )_{q} = \left ( {n \atop i}
\right )_{q} + q^{n+1-i} \left ( {n \atop i-1} \right )$
enables to prove
 the equality $ \left ( {n \atop i} \right )_{q} =
\frac{n_{q} !}{(n-i)_{q}! i_{q} !}$ where $n_{q} =
1+q+\cdots+q^{n-1}$ and $n_{q}! = n_{q}(n-1)_{q} \cdots
1_{q}$ (see for instance \cite{ka}).

\begin{Pro} Let $E^n_{i}$ denotes the path of
length $n$ with source vertex $K^i$ (the index $i$ is to be
considered modulo the order of $G$ in case $G$ is finite).
Then $$ E^n_{i}. E^m_{j} = q^{jn} \left ( {n+m \atop n}
\right )_{q} E^{n+m}_{i+j}. $$
\end{Pro}

\noindent \textbf{Proof} : We consider the product
$E^n_{0}.E^m_{0}$. Using the formula of the Theorem this
product is a multiple of $E^{n+m}_{0}$ by a scalar which we
denote $X^{n+m}_{n}$. By induction $X^{n+m}_{n} = \left (
{n+m \atop n} \right )_{q}$ since $X_{1}^{m+1} = \left (
{m+1 \atop 1} \right )_{q}$ and the associativity formula
$$ E^{1}_{0}.(E^{n}_{0}.E^m_{0}) =
(E^{1}_{0}.E^n_{0}).E^m_{0} $$ provides the recursive
formula $ \left ( {n+m+1 \atop 1} \right )_{q} X^{n+m}_{n}
= \left ( {n+1 \atop 1} \right )_{q} X_{n+1}^{n+m+1} $
which is clearly satisfied while replacing $X_{n}^{n+m}$ by
$ \left ( {n+m \atop n} \right )_{q}$. Finally we record
that $K^i.E^n_{0} = E^n_{i}$ and $E^n_{0}.K^j =
q^{nj}E^n_{j}$. Hence $$ E^n_{i} . E^m_{j} = K^i.E^n_{0} .
K^j.E^m_{0} = q^{nj} K^{i+j}.E^n_{0}.E^m_{0}. $$

\section{\sf Hopf bimodules}

The proof of the structure Theorem of the preceding section
makes use of the classification of Hopf bimodules over a
group algebra obtained in \cite{ciro}, see also
\cite{ci,angr,diparo}. Let $k$ be a field, $G$ be a group,
$\mathcal{C}$ be the set of conjugacy classes, $Z_{C}$ be
the centralizer of some element in the conjugacy $C$ and
mod $kZ_{C}$ be the category of right $kZ_{C}$-modules.

\begin{Theo}
The category of Hopf bimodules $\mathcal{B}(kG)$ is
equivalent to the cartesian product of categories $\prod_{C
\in \mathcal{C}} \ {\rm mod}\ kZ_{C}$.
\end{Theo}

For a complete proof we refer to the quoted references. In
order to comply with the requirements of the preceding
section, we describe the functor $W$ which associates to a
Hopf bimodule $B$ the family $\{^{u(C)}B^{1}\}_{C\in
\mathcal{C}}$ where $u(C)$ is some element in $C$ and $$
^{u(C)}\!B^{1} = \{b \in B | \delta_{L}(b) =u(C) \otimes b
\quad {\rm and}\quad \delta_{R}(b) = b \otimes 1 \}. $$
Note that $^{u(C)}B^{1}$ is a right $kZ_{C}$-module. In
this setting the corresponding ramification data is
 $$ r = \sum_{C
\in\mathcal{C}} (\dim_{k} ^{u(C)}\!\!B^{1} )C.$$
 In order to provide
at least one Hopf algebra structure to the path coalgebra of a
Hopf quiver with ramification data $r$, one can chose a family
$\{M_{C}\}_{C \in \mathcal{C}}$ of trivial $kZ_{C}$-modules of
adapted dimension, namely $ \dim_{k} M_{C} = r_{C} $. More
generally, the following result is now evident :

\begin{Theo}
Let $G$ be a group, $r$ be a ramification data and $Q$ be the
corresponding Hopf quiver. The complete list of graded Hopf
algebra structures on the path coalgebra $kQ$ is in one to one
correspondence with the set of collections $\{M_{C}\}_{C \in
\mathcal{C}}$ where $M_{C}$ is a right $kZ_{u(C)}$-module of
dimension $r_{C}$.
\end{Theo}

\begin{Ex}\sf
Consider $G$ a cyclic group (finite or not) with a generator $K$,
and the ramification data $r = K$. In other words, only one
conjugacy class is highlighted by $r$ and have coefficient one.
Hence we consider the set of one-dimensional $kG$-modules, in
order to obtain the complete list of Hopf algebras structures on
the path coalgebra of the corresponding Hopf quiver.

A one-dimensional $kG$-module is a non-zero scalar $q$ if
$G$ is infinite and a $|G|$-root of unity $q$ if $G$ is
finite. The resulting Hopf bimodule is obtained through the
functor $V$ of \cite{ci}, and is the Hopf bimodule we have
described in the last example of the preceding section.
\end{Ex}

\footnotesize \noindent C.C.:\\ Universit\'e de Montpellier
2, D\'epartement de Math\'ematiques, \\F-34095 Montpellier
cedex 5, France \\{\tt Claude.Cibils@math.univ-montp2.fr}

\vskip3mm \noindent M.R.: \\Ecole Normale
Sup\'erieure\\D\'epartement de math\'ematiques et
applications\\F-75230 Paris cedex 05, France.
\\ {\tt Marc.Rosso@ens.fr}

\noindent September 2000


\begin{thebibliography}{99}

\footnotesize


\bibitem{angr}
Andruskiewitsch, N.; Gra\~na, M. Braided Hopf algebras over
non-abelian finite groups. Colloquium on Operator Algebras
and Quantum Groups (Vaquer\'\i as, 1997). Bol. Acad. Nac.
Cienc. (C\'ordoba) 63 (1999), 45--78.

\bibitem{auresm} Auslander, M.; Reiten, I.; Smal\o, S.
Representation theory of Artin algebras. Cambridge Studies in
Advanced Mathematics. 36. Cambridge: Cambridge University Press
(1995).


\bibitem{ci}
Cibils, C. Tensor product of Hopf bimodules over a group.
Proc. Amer. Math. Soc. 125 (1997), 1315--1321.

\bibitem{cimsri}
Cibils, C. Hopf quivers. Lecture on video at MSRI (1999).\\
{\scriptsize \tt
http://www.msri.org/publications/ln/msri/1999/hopfalg/cibils/1/index.html}

\bibitem{ciro}
Cibils, C.; Rosso, M. Alg\`{e}bres des chemins quantiques.
Adv. Math. 125 (1997), 171--199.

\bibitem{chmo}
Chin, W.; Montgomery, S. Basic coalgebras. Modular interfaces
(Riverside, CA, 1995), 41--47, AMS/IP Stud. Adv. Math., 4, Amer.
Math. Soc., Providence, RI, 1997.

\bibitem{diparo}
Dijkgraaf, R.; Pasquier, V.; Roche, P. Quasi Hopf algebras,
group cohomology and orbifold models. Recent advances in
field theory (Annecy-le-Vieux, 1990). Nuclear Phys. B Proc.
Suppl. 18B (1990), 60--72 (1991)

\bibitem{gr}
Green, J. A. Locally finite representations. J. Algebra 41
(1976), 137--171

\bibitem{grso}
Green, E. L.; Solberg, \O .
 Basic Hopf algebras and quantum groups. Math. Z. 229 (1998), 45--76

\bibitem{ka}
Kassel, C. Quantum groups. Graduate Texts in Mathematics, 155.
Springer-Verlag, New York, 1995

\bibitem{mo}
Montgomery, S. Hopf algebras and their actions on rings.
CBMS Regional Conference Series in Mathematics, 82.
Published for the Conference Board of the Mathematical
Sciences, Washington, DC; by the American Mathematical
Society, Providence, RI, 1993

\bibitem{ni} Nichols W.D. Pointed irreducible algebras. J. of
Algebra 57 (1979), 64--76

\bibitem{ra} Radford, D. E. Finite-dimensional simple-pointed Hopf
algebras. J. Algebra 211 (1999), 686--710

\bibitem{ro} Rosso, Marc Quantum groups and quantum
shuffles. Invent. Math. 133 (1998), 399--416

\bibitem{sw}
Sweedler, M.E.: Hopf algebras. Benjamin. New York U.S.A.
1969

\bibitem{ta}
Takeuchi, M. Free Hopf algebras generated by coalgebras. J. Math.
Soc. Japan 23 (1971), 561--582.

\bibitem{tooyzh} Torrecillas, B.; Van Oystaeyen, F.; Zhang, Y. H.
 Coflat monomorphisms of coalgebras. J. Pure Appl. Algebra 128 (1998), 171--183.

\end{thebibliography}
\end{document}